# NONCONCAVE PENALIZED LIKELIHOOD WITH A DIVERGING NUMBER OF PARAMETERS

By Jianqing Fan[1] and Heng Peng

*Princeton University and The Chinese University of Hong Kong*

A class of variable selection procedures for parametric models via nonconcave penalized likelihood was proposed by Fan and Li to simultaneously estimate parameters and select important variables. They demonstrated that this class of procedures has an oracle property when the number of parameters is finite. However, in most model selection problems the number of parameters should be large and grow with the sample size. In this paper some asymptotic properties of the nonconcave penalized likelihood are established for situations in which the number of parameters tends to $\infty$ as the sample size increases. Under regularity conditions we have established an oracle property and the asymptotic normality of the penalized likelihood estimators. Furthermore, the consistency of the sandwich formula of the covariance matrix is demonstrated. Nonconcave penalized likelihood ratio statistics are discussed, and their asymptotic distributions under the null hypothesis are obtained by imposing some mild conditions on the penalty functions. The asymptotic results are augmented by a simulation study, and the newly developed methodology is illustrated by an analysis of a court case on the sexual discrimination of salary.

## 1. Introduction.

1.1. *Background.* The idea of penalization is very useful in statistical modeling, particularly in variable selection, which is fundamental to the field. Most traditional variable selection procedures, such as Akaike's information criterion AIC [Akaike (1973)], Mallows' $C_p$ [Mallows (1973)] and the Bayesian information criterion BIC [Schwarz (1978)], use a fixed penalty on

Received October 2002; revised April 2003.
[1]Supported by NSF Grant DMS-02-04329 and NIH Grant R01 HL69720.

the size of a model. Some new variable selection procedures suggest the use of a data adaptive penalty to replace fixed penalties [i.e., Bai, Rao and Wu (1999) and Shen and Ye (2002)]. However, all these procedures follow stepwise and subset selection procedures to select variables. Stepwise and subset selection procedures make these procedures computationally intensive, hard to derive sampling properties, and unstable [see, e.g., Breiman (1996) and Fan and Li (2001)]. In contrast, most convex penalties, such as quadratic penalties, often produce shrinkage estimators of parameters that make tradeoffs between bias and variance such as those in smoothing splines. However, they can create unnecessary biases when the true parameters are large and parsimonious models cannot be produced.

To overcome the inefficiency of traditional variable selection procedures, Fan and Li (2001) proposed a unified approach via nonconcave penalized least squares to automatically and simultaneously select variables and estimate the coefficients of variables. This method not only retains the good features of both subset selection and ridge regression, but also produces sparse solutions (many estimated coefficients are 0), ensures continuity of the selected models (for the stability of model selection) and has unbiased estimates for large coefficients. This is achieved by choosing suitable penalized nonconcave functions, such as the smoothly clipped absolute deviation (SCAD) penalty that was proposed by Fan (1997) (to be defined in Section 2). Other penalized least squares, such as the bridge regression proposed by Frank and Friedman (1993) and Lasso proposed by Tibshirani (1996, 1997), can also be studied under this unified work. The nonconcave penalized least-squares approach also corresponds to a Bayesian model selection with an improper prior and can be easily extended to likelihood-based models in various statistical contexts, such as generalized linear modeling, nonparametric modeling and survival analysis. For example, Antoniadis and Fan (2001) used this approach in wavelet analysis, and Fan and Li (2002) applied the technique to the Cox proportional hazards model and the frailty model.

1.2. *Nonconcave penalized likelihood.* One distinguishing feature of the nonconcave penalized likelihood approach is that it can simultaneously select variables and estimate coefficients of variables. This enables us to establish the sampling properties of the nonconcave penalized likelihood estimates.

Let $\log f(V, \beta)$ be the underlying likelihood for a random vector $V$. This includes the likelihood of the form $\ell(X^T \beta, Y)$ of the generalized linear model [McCullagh and Nelder (1989)]. Let $p_\lambda(|\beta_j|)$ be a nonconcave penalized function that is indexed by a regularization parameter $\lambda$. The penalized likelihood estimator then maximizes

$$\sum_{i=1}^{n} \log f(V_i, \beta) - \sum_{j=1}^{p} p_\lambda(|\beta_j|). \tag{1.1}$$

The parameter $\lambda$ can be chosen by cross-validation [see Breiman (1996) and Tibshirani (1996)].

Various algorithms have been proposed to optimize such a high-dimensional nonconcave likelihood function. The modified Newton–Raphson algorithm was proposed by Fan and Li (2001). The idea of the graduated nonconvexity algorithm was proposed by Blake and Zisserman (1987) and Blake (1989), and was implemented by Nikolova, Idier and Mohammad-Djafari (1998). Tibshirani (1996, 1997) and Fu (1998) proposed different algorithms for the $L_p$-penalty. One can also use a stochastic optimization method, such as simulated annealing. See Geman and Geman (1984) and Gilks, Richardson and Spiegelhalter (1996) for more discussions.

For the finite parameter case, Fan and Li (2001) established an "oracle property," to use the terminology of Donoho and Johnstone (1994). If there were an oracle assisting us in selecting variables, then we would select variables only with nonzero coefficients and apply the MLE to this submodel and estimate the remaining coefficients as 0. This ideal estimator is called an oracle estimator. Fan and Li (2001) demonstrated that penalized likelihood estimators are asymptotically as efficient as this ideal oracle estimator for certain penalty functions, such as SCAD and the hard thresholding penalty. Fan and Li (2001) also proposed a sandwich formula for estimating the standard error of the estimated nonzero coefficients and empirically verifying the consistency of the formula. Knight and Fu (2000) studied the asymptotic behavior of the Lasso type of estimator. Under some appropriate conditions, they showed that the limiting distributions have positive probability mass at 0 when the true value of the parameters is 0, and they established asymptotic normality for large parameters in some sense.

1.3. *Real issues in model selection.* In practice, many variables are introduced to reduce possible modeling biases. The number of introduced variables depends on the sample size, which reflects the estimability of the parametric problem.

An early reference on this kind of problem is the seminal paper of Neyman and Scott (1948). In the early years, from problems in X-ray crystallography, where the typical values for the number of parameters $p$ and sample size $n$ are in the ranges 10 to 500 and 100 to 10,000, respectively, Huber (1973) noted that in a variable selection context the number of parameters is often large and should be modeled as $p_n$, which tends to $\infty$. Now, with the advancement of technology and huge investment in various forms of data gathering, as Donoho (2000) demonstrated with web term-document data, gene expression data and consumer financial history data, large sample sizes with high dimensions are important characteristics. He also observed that even in a classical setting such as the Framingham heart study, the sample size is as large as $N = 25{,}000$ and the dimension is $p = 100$, which can be modeled as $p = O(n^{1/2})$ or $p = O(n^{1/3})$.

Nonparametric regression is another class of examples that uses diverging parameters. In spline modeling an unknown function is frequently approximated by its finite series expansion with the number of parameters depending on the sample size. In regression splines, Stone, Hansen, Kooperberg and Truong (1997) regard nonparametric problems as large parametric problems and extend traditional variable selection techniques to select important terms. Smoothing splines can also be regarded as a large parametric problem [Green and Silverman (1994)]. To achieve the stability of the resulting estimate (e.g., smoothness), instead of selecting variables a quadratic penalty is frequently used to shrink the estimated parameters [Cox and O'Sullivan (1990)]. Thus, our formulation and results have applications to the problem of nonparametric estimation.

Fan and Li (2001) laid down important groundwork on variable selection problems, but their theoretical results are limited to the finite-parameter setting. While their results are encouraging, the fundamental problems with a growing number of parameters have not been addressed. In fact, the full advantages of the penalized likelihood method in model selection have not been convincingly demonstrated. For example, for finite-parameter problems, owing to the root-$n$-consistency of estimated parameters, many naive and simple model selection procedures also possess the oracle property. To wit, a simple thresholding estimator such as $\hat{\beta}_j I(|\hat{\beta}_j| > n^{-1/4})$, which completely ignores the correlation structure and the scale of the parameter, also possesses the oracle property. Thus, it is uncertain whether the oracle property of Fan and Li (2001) is genuine to the penalized likelihood method or an artifact of the finite-parameter formulation.

To this end, we consider the log-likelihood series $\log f_n(V_n, \beta_n)$, where $f_n(V_n, \beta_n)$ is the density of the random variable $V_n$, all of which relate to the sample size $n$, and assume without loss of generality that, unknown to us, the first $s_n$ components of $\beta_n$, denoted by $\beta_{n1}$, do not vanish and the remaining $p_n - s_n$ coefficients, denoted by $\beta_{n2}$, are 0. Our objectives in this paper are to investigate the following asymptotic properties of a nonconcave penalized likelihood estimator.

1. (Oracle property.) Under certain conditions of the likelihood function and for certain penalty functions (e.g., SCAD), if $p_n$ does not grow too fast, then by the proper choice of $\lambda_n$ there exists a penalized likelihood estimator such that $\hat{\beta}_{n2} = 0$ and $\hat{\beta}_{n1}$ behaves the same as the case in which $\beta_{n2} = 0$ is known in advance.
2. (Asymptotic normality.) As the length of $\hat{\beta}_{n1}$ depends on $n$, we will consider its arbitrary linear combination $A_n \hat{\beta}_{n1}$, where $A_n$ is a $q \times s_n$ matrix for any finite $q$. We will show that this linear combination is asymptotically normal. Furthermore, let $\hat{\beta}_{n1}^o$ be the oracle estimator, thus maximizing the likelihood of the ideal submodel $\sum_{i=1}^n \log f_n(V_i, \beta_{n1})$. We will show that $A_n \hat{\beta}_{n1}^o$ is also asymptotically normal. We will study the conditions under which the two covariance matrices are identical. This will demonstrate the oracle property mentioned above.

3. (Consistency of the sandwich formula.) Let $\hat{\Sigma}_n$ be an estimated covariance matrix for $\hat{\beta}_{n1}$, using the sandwich formula based on the penalized likelihood (1.1). We will show that the covariance matrix $\hat{\Sigma}_n$ is a consistent estimate in the sense that $A_n^T \hat{\Sigma}_n A_n$ converges to the $q \times q$ asymptotic covariance matrix of $A_n \hat{\beta}_{n1}$.
4. (Likelihood ratio theory.) If one tests the linear hypothesis $H_0 : A_n \hat{\beta}_{n1} = 0$ and uses the twice-penalized likelihood ratio statistic, then this statistic asymptotically follows a $\chi^2$ distribution.

The asymptotic properties of any finite components of $\hat{\beta}$ are included in the above formulation by taking a special matrix $A_n$. Furthermore, the asymptotic properties and variable selection of linear components in any partial linear model can be analyzed this way if we use a series such as a Fourier series or polynomial splines to estimate the nonparametric component.

1.4. *Outline of the paper.* In Section 2 we briefly review the nonconcave penalized likelihood. The asymptotic results of penalized likelihood are presented in Section 3. We discuss the conditions that are imposed on the likelihood and penalty functions in Section 3.1 and present our main results in Sections 3.2–3.4. An application of the proposed methodology and a simulation study are presented in Section 4. The proofs of our results are given in Section 5. Technical details are relegated to the Appendix.

**2. Penalty function.** Penalty functions largely determine the sampling properties of the penalized likelihood estimators. To select a good penalty function, Fan and Li (2001) proposed three principles that a good penalty function should satisfy: unbiasedness, in which there is no overpenalization of large parameters to avoid unnecessary modeling biases; sparsity, as the resulting penalized likelihood estimators should follow a thresholding rule such that insignificant parameters are automatically set to 0 to reduce model complexity; and continuity to avoid instability in model prediction, whereby the penalty function should be chosen such that its corresponding penalized likelihood produces continuous estimators of data. More details can be found in the work of Fan and Li (2001) and Antoniadis and Fan (2001).

To gain some insight into the choice of penalty functions, let us first consider a simple form of (1.1), that is, the penalized least-squares problem:

$$\tfrac{1}{2}(z-\theta)^2 + p_\lambda(|\theta|).$$

It is well known that the $L_2$-penalty $p_\lambda(|\theta|) = \lambda|\theta|^2$ leads to a ridge regression. A generalization is the $L_q$-penalty $p_\lambda(|\theta|) = \lambda|\theta|^q$, $q > 1$. These penalties reduce variability via shrinking the solutions, but do not have the properties of sparsity.

The $L_1$-penalty $p_\lambda(|\theta|) = \lambda|\theta|$ yields a soft thresholding rule

$$\hat{\theta} = \operatorname{sgn}(z)(|z| - \lambda)_+.$$

Tibshirani (1996, 1997) applied the $L_1$-penalty to a general least-squares and likelihood setting. Knight and Fu (2000) studied the $L_q$-penalty when $q < 1$. While the $L_q$-penalty ($q \leq 1$) functions result in sparse solutions, they cannot keep the resulting estimators unbiased for large parameters due to excessive penalty at large values of parameters. Another type of penalty function is the hard thresholding penalty function

$$p_\lambda(|\theta|) = \lambda^2 - (|\theta| - \lambda)^2 I(|\theta| < \lambda),$$

which results in the hard thresholding rule [see Antoniadis (1997) and Fan (1997)]

$$\hat{\theta} = zI(|z| > \lambda),$$

but the estimator is not continuous in the data $z$.

As the penalty functions above cannot simultaneously satisfy the aforementioned three principles, motivated by wavelet analysis, Fan (1997) proposed a continuous differentiable penalty function called the smoothly clipped absolute deviation (SCAD) penalty, which is defined by

$$p'_\lambda(\theta) = \lambda \left\{ I(\theta \leq \lambda) + \frac{(a\lambda - \theta)_+}{(a-1)\lambda} I(\theta > \lambda) \right\} \qquad \text{for some } a > 2 \text{ and } \theta > 0.$$

The solution for this penalty function is given by

$$\hat{\theta} = \begin{cases} \operatorname{sgn}(z)(|z| - \lambda)_+, & \text{when } |z| \leq 2\lambda, \\ \{(a-1)z - \operatorname{sgn}(z)a\lambda\}/(a-2), & \text{when } 2\lambda \leq |z| \leq a\lambda, \\ z, & \text{when } |z| > a\lambda. \end{cases}$$

The solution satisfies the three properties that were proposed by Fan and Li (2001).

**3. Properties of penalized likelihood estimation.** In this section we study the sampling properties of the penalized likelihood estimators proposed in Section 1 in the situation where the number of parameters tends to $\infty$ with increasing sample size. We discuss some conditions of the penalty and likelihood functions in Section 3.1 and show their differences from those under finite parameters. Though the imposed conditions are not the weakest possible, they make technical analysis easily understandable. Our main results are presented in Section 3.2.

### 3.1. *Regularity conditions.*

3.1.1. *Regularity condition on penalty.* Let $a_n = \max_{1 \leq j \leq p_n} \{p'_{\lambda_n}(|\beta_{n0j}|), \beta_{n0j} \neq 0\}$ and $b_n = \max_{1 \leq j \leq p_n} \{p''_{\lambda_n}(|\beta_{n0j}|), \beta_{n0j} \neq 0\}$. Then we need to place the following conditions on the penalty functions:

(A) $\liminf_{n \to +\infty} \liminf_{\theta \to 0+} p'_{\lambda_n}(\theta)/\lambda_n > 0$;
(B) $a_n = O(n^{-1/2})$;
(B′) $a_n = o(1/\sqrt{np_n})$;
(C) $b_n \to 0$ as $n \to +\infty$;
(C′) $b_n = o_p(1/\sqrt{p_n})$;
(D) there are constants $C$ and $D$ such that, when $\theta_1, \theta_2 > C\lambda_n$, $|p''_{\lambda_n}(\theta_1) - p''_{\lambda_n}(\theta_2)| \leq D|\theta_1 - \theta_2|$.

Condition (A) makes the penalty function singular at the origin so that the penalized likelihood estimators possess the sparsity property. Conditions (B) and (B′) ensure the unbiasedness property for large parameters and the existence of the root-$n$-consistent penalized likelihood estimator. Conditions (C) and (C′) guarantee that the penalty function does not have much more influence than the likelihood function on the penalized likelihood estimators. Condition (D) is a smoothness condition that is imposed on the nonconcave penalty functions. Under the condition (H) all of these conditions are satisfied by the SCAD penalty and the hard thresholding penalty, as $a_n = 0$ and $b_n = 0$ when $n$ is large enough.

3.1.2. *Regularity conditions on likelihood functions.* Due to the diverging number of parameters, we cannot assume that likelihood functions are invariant in our study. Some conditions have to be strengthened to keep uniform properties for the likelihood functions and sample series. A higher-order moment of the likelihood functions is a simple and direct method to keep uniform properties, as compared to the usual conditions in the asymptotic theory of the likelihood estimate under finite parameters [see, e.g., Lehmann (1983)]. The conditions that are imposed on the likelihood functions are as follows:

(E) For every $n$ the observations $\{V_{ni}, i = 1, 2, \ldots, n\}$ are independent and identically distributed with the probability density $f_n(V_{n1}, \beta_n)$, which has a common support, and the model is identifiable. Furthermore, the first and second derivatives of the likelihood function satisfy the equations

$$E_{\beta_n}\left\{\frac{\partial \log f_n(V_{n1}, \beta_n)}{\partial \beta_{nj}}\right\} = 0 \qquad \text{for } j = 1, 2, \ldots, p_n$$

and

$$E_{\beta_n}\left\{\frac{\partial \log f_n(V_{n1},\beta_n)}{\partial \beta_{nj}}\frac{\partial \log f_n(V_{n1},\beta_n)}{\partial \beta_{nk}}\right\} = -E_{\beta_n}\left\{\frac{\partial^2 \log f_n(V_{n1},\beta_n)}{\partial \beta_{nj}\,\partial \beta_{nk}}\right\}.$$

(F) The Fisher information matrix

$$I_n(\beta_n) = E\left[\left\{\frac{\partial \log f_n(V_{n1},\beta_n)}{\partial \beta_n}\right\}\left\{\frac{\partial \log f_n(V_{n1},\beta_n)}{\partial \beta_n}\right\}^T\right]$$

satisfies conditions

$$0 < C_1 < \lambda_{\min}\{I_n(\beta_n)\} \le \lambda_{\max}\{I_n(\beta_n)\} < C_2 < \infty \qquad \text{for all } n$$

and, for $j,k = 1,2,\ldots,p_n$,

$$E_{\beta_n}\left\{\frac{\partial \log f_n(V_{n1},\beta_n)}{\partial \beta_{nj}}\frac{\partial \log f_n(V_{n1},\beta_n)}{\partial \beta_{nk}}\right\}^2 < C_3 < \infty$$

and

$$E_{\beta_n}\left\{\frac{\partial^2 \log f_n(V_{n1},\beta_n)}{\partial \beta_{nj}\,\partial \beta_{nk}}\right\}^2 < C_4 < \infty.$$

(G) There is a large enough open subset $\omega_n$ of $\Omega_n \in R^{p_n}$ which contains the true parameter point $\beta_n$, such that for almost all $V_{ni}$ the density admits all third derivatives $\partial f_n(V_{ni},\beta_n)/\partial \beta_{nj}\beta_{nk}\beta_{nl}$ for all $\beta_n \in \omega_n$. Furthermore, there are functions $M_{njkl}$ such that

$$\left|\frac{\partial \log f_n(V_{ni},\beta_n)}{\partial \beta_{nj}\beta_{nk}\beta_{nl}}\right| \le M_{njkl}(V_{ni})$$

for all $\beta_n \in \omega_n$, and

$$E_{\beta_n}\{M_{njkl}^2(V_{ni})\} < C_5 < \infty$$

for all $p,n$ and $j,k,l$.

(H) Let the values of $\beta_{n01},\beta_{n02},\ldots,\beta_{n0s_n}$ be nonzero and $\beta_{n0(s_n+1)},\beta_{n02},\ldots,\beta_{n0p_n}$ be zero. Then $\beta_{n01},\beta_{n02},\ldots,\beta_{n0s_n}$ satisfy

$$\min_{1 \le j \le s_n} |\beta_{n0j}|/\lambda_n \to \infty \qquad \text{as } n \to \infty.$$

Under conditions (F) and (G), the second and fourth moments of the likelihood function are imposed. The information matrix of the likelihood function is assumed to be positive definite, and its eigenvalues are uniformly bounded. These conditions are stronger than those of the usual asymptotic likelihood theory, but they facilitate the technical derivations.

Condition (H) seems artificial, but it is necessary for obtaining the oracle property. In a finite-parameter situation this condition is implicitly assumed, and is in fact stronger than that imposed here. Condition (H) explicitly shows the rate at which the penalized likelihood can distinguish nonvanishing parameters from 0. Its zero component can be relaxed as

$$\max_{s_n+1\leq j\leq p_n} |\beta_{n0j}|/\lambda_n \to 0 \qquad \text{as } n \to \infty.$$

3.2. *Oracle properties.* Recall that $V_{ni}, i=1,\ldots,n$, are independent and identically distributed random variables with density $f_n(V_n, \beta_{n0})$. Let

$$L_n(\beta_n) = \sum_{i=1}^n \log f_n(V_{ni}, \beta_n)$$

be the log-likelihood function and let

$$Q_n(\beta_n) = L_n(\beta_n) - n\sum_{j=1}^{p_n} p_{\lambda_n}(|\beta_{nj}|)$$

be the penalized likelihood function.

THEOREM 1 (Existence of penalized likelihood estimator). *Suppose that the density $f_n(V_n, \beta_{n0})$ satisfies conditions* (E)–(G), *and the penalty function $p_{\lambda_n}(\cdot)$ satisfies conditions* (B)–(D). *If $p_n^4/n \to 0$ as $n \to \infty$, then there is a local maximizer $\hat{\beta}_n$ of $Q(\beta_n)$ such that $\|\hat{\beta}_n - \beta_{n0}\| = O_p\{\sqrt{p_n}(n^{-1/2} + a_n)\}$, where $a_n$ is given in Section* 3.1.1.

It is easy to see that if $a_n$ satisfies condition (B), that is, $a_n = O(n^{-1/2})$, then there is a root-$(n/p_n)$-consistent estimator. This consistent rate is the same as the result of the M-estimator that was studied by Huber (1973), in which the number of parameters diverges. The convergence rate of $a_n$ for the usual convex penalties, such as the $L_q$-penalty with $q \geq 1$, largely depends on the convergence rate of $\lambda_n$. As these penalties do not have an unbiasedness property, they require that $\lambda_n$ satisfy the condition $\lambda_n = O(n^{-1/2})$ in order to have a root-$(n/p_n)$-consistent estimator for the penalized likelihood estimator. This requirement will make it difficult to choose $\lambda_n$ for penalized likelihood in practice. However, if the penalty function is a SCAD or hard thresholding penalty, and condition (H) is satisfied by the model, it is clear that $a_n = 0$ when $n$ is large enough. The root-$(n/p_n)$-consistent penalized likelihood estimator indeed exists with probability tending to 1, and no requirements are imposed on the convergence rate of $\lambda_n$.

Denote

$$\Sigma_{\lambda_n} = \text{diag}\{p''_{\lambda_n}(\beta_{n01}),\ldots,p''_{\lambda_n}(\beta_{n0s_n})\}$$

and

$$\mathbf{b}_n = \{p'_{\lambda_n}(|\beta_{n01}|)\,\text{sgn}(\beta_{n01}),\ldots,p'_{\lambda_n}(|\beta_{n0s_n}|)\,\text{sgn}(\beta_{n0s_n})\}^T.$$

THEOREM 2 (Oracle property). *Under conditions* (A)–(H), *if* $\lambda_n \to 0$, $\sqrt{n/p_n}\lambda_n \to \infty$ *and* $p_n^5/n \to 0$ *as* $n \to \infty$, *then, with probability tending to* 1, *the root-$(n/p_n)$-consistent local maximizer* $\hat{\beta}_n = \binom{\hat{\beta}_{n1}}{\hat{\beta}_{n2}}$ *in Theorem* 1 *must satisfy:*

(i) (Sparsity) $\hat{\beta}_{n2} = 0$.
(ii) (Asymptotic normality)

$$\sqrt{n} A_n I_n^{-1/2}(\beta_{n01})\{I_n(\beta_{n01}) + \Sigma_{\lambda_n}\}$$
$$\times [\hat{\beta}_{n1} - \beta_{n01} + \{I_n(\beta_{n01}) + \Sigma_{\lambda_n}\}^{-1}\mathbf{b}_n] \xrightarrow{\mathcal{D}} \mathcal{N}(0, G),$$

where $A_n$ is a $q \times s_n$ matrix such that $A_n A_n^T \to G$, and $G$ is a $q \times q$ nonegative symmetric matrix.

By Theorem 2 the sparsity and the asymptotic normality are still valid when the number of parameters diverges. In fact, the oracle property holds for the SCAD and the hard thresholding penalty function. When $n$ is large enough, $\Sigma_{\lambda_n} = 0$ and $\mathbf{b}_n = 0$ for the SCAD and the hard thresholding penalty. Hence, Theorem 2(ii) becomes

$$\sqrt{n} A_n I_n^{1/2}(\beta_{n01})(\hat{\beta}_{n1} - \beta_{n01}) \xrightarrow{\mathcal{D}} \mathcal{N}(0, G),$$

which has the same efficiency as the maximum likelihood estimator of $\beta_{n01}$ based on the submodel with $\beta_{n02} = 0$ known in advance. This demonstrates that the nonconcave penalized likelihood estimator is as efficient as the oracle one. Intrinsically, unbiasedness and singularity at the origin of the SCAD and the hard thresholding penalty functions guarantee this sampling property.

The $L_q$-penalty, $q \geq 1$, cannot simultaneously satisfy the conditions $\lambda_n = O_p(n^{-1/2})$ and $\sqrt{n/p}\lambda_n \to \infty$ as $n \to \infty$. These penalty functions cannot produce estimators with the oracle property. The $L_q$-penalty, $q < 1$, may satisfy these two conditions at same time. As shown by Knight and Fu (2000) in a finite-parameter setting, it might also have sampling properties that are similar to the oracle property when the number of parameters diverges. However, the bias term in Theorem 2(ii) cannot be ignored.

The condition $p_n^4/n \to 0$ or $p_n^5/n \to 0$ as $n \to \infty$ seems somewhat strong. By refining the structure of the log-likelihood function, such as the generalized linear model $\ell(X^T\beta, Y)$ or the $M$-estimator from $\sum_{i=1}^n \rho(Y_i - X_i^T\beta)$, the condition can be weakened to $p_n^3/n \to 0$ as $n \to \infty$. This condition is in line with that of Huber (1973).

3.3. *Estimation of covariance matrix.* As in Fan and Li (2001), by the sandwich formula let

$$\hat{\Sigma}_n = n\{\nabla^2 L_n(\hat{\beta}_{n1}) - n\Sigma_{\lambda_n}(\hat{\beta}_{n1})\}^{-1}$$
$$\times \widehat{\mathrm{cov}}\{\nabla L_n(\hat{\beta}_{n1})\}\{\nabla^2 L_n(\hat{\beta}_{n1}) - n\Sigma_{\lambda_n}(\hat{\beta}_{n1})\}^{-1}$$

be the estimated covariance matrix of $\hat{\beta}_{n1}$, where

$$\widehat{\mathrm{cov}}\{\nabla L_n(\hat{\beta}_{n1})\} = \left\{\frac{1}{n}\sum_{i=1}^n \frac{\partial L_{ni}(\hat{\beta}_{n1})}{\partial \beta_j}\frac{\partial L_{ni}(\hat{\beta}_{n1})}{\partial \beta_k}\right\}$$
$$- \left\{\frac{1}{n}\sum_{i=1}^n \frac{\partial L_{ni}(\hat{\beta}_{n1})}{\partial \beta_j}\right\}\left\{\frac{1}{n}\sum_{i=1}^n \frac{\partial L_{ni}(\hat{\beta}_{n1})}{\partial \beta_k}\right\}.$$

Denote by

$$\Sigma_n = \{I_n(\beta_{n01}) + \Sigma_{\lambda_n}(\beta_{n01})\}^{-1} I_n(\beta_{n01}) \{I_n(\beta_{n01}) + \Sigma_{\lambda_n}(\beta_{n01})\}^{-1}$$

the asymptotic variance of $\hat{\beta}_{n1}$ in Theorem 2(ii).

THEOREM 3 (Consistency of the sandwich formula). *If conditions* (A)–(H) *are satisfied and* $p_n^5/n \to 0$ *as* $n \to \infty$, *then we have*

$$(3.1) \qquad A_n \hat{\Sigma}_n A_n^T - A_n \Sigma_n A_n^T \overset{p}{\to} 0 \qquad as\ n \to \infty$$

*for any* $q \times s_n$ *matrix* $A_n$ *such that* $A_n A_n^T = G$, *where* $q$ *is any fixed integer.*

Theorem 3 not only proves a conjecture of Fan and Li (2001) about the consistency of the sandwich formula for the standard error matrix, but also extends the result to the situation with a growing number of parameters. The consistent result also offers a way to construct a confidence interval for the estimates of parameters. For a review of sandwich covariance matrix estimation, see the paper of Kauermann and Carroll (2001).

3.4. *Likelihood ratio test.* One of the most celebrated methods in statistics is the likelihood ratio test. Can it also be applied to the penalized likelihood context with a diverging number of parameters? To answer this question, consider the problem of testing linear hypotheses:

$$H_0 : A_n \beta_{n01} = 0 \quad \text{vs.} \quad H_1 : A_n \beta_{n01} \neq 0,$$

where $A_n$ is a $q \times s_n$ matrix and $A_n A_n^T = I_q$ for a fixed $q$. This problem includes testing simultaneously the significance of a few covariate variables.

In the penalized likelihood context a natural likelihood ratio test for the problem is

$$T_n = 2\left\{\sup_{\Omega_n} Q(\beta_n \mid \mathbf{V}) - \sup_{\Omega_n, A_n \beta_{n1} = 0} Q(\beta_n \mid \mathbf{V})\right\}.$$

The following theorem drives the asymptotic null distribution of the test statistic. It shows that the classical likelihood theory continues to hold for the problem with a growing number of parameters in the penalized likelihood context. It enables one to apply the traditional likelihood ratio method for testing linear hypotheses. In particular, it allows one to simultaneously test whether a few variables are statistically significant by taking some specific matrix $A_n$.

THEOREM 4. *When conditions* (A)–(H), $(B')$ *and* $(C')$ *are satisfied, under* $H_0$ *we have*

$$T_n \xrightarrow{\mathcal{D}} \chi_q^2, \tag{3.2}$$

*provided that* $p_n^5/n \to 0$ *as* $n \to \infty$.

For the usual likelihood without penalization, Portnoy (1988) and Murphy (1993) showed that the Wilks type of result continues to hold for specific problems. Our results can be regarded as a further generalization of theirs.

**4. Numerical examples.** In this section we illustrate the techniques of our method via an analysis of a data set in a lawsuit and verify the finite-sample performance via a simulation experiment.

4.1. *A real data example.* The Fifth National Bank of Springfield faced a gender discrimination suit. [This example and the accompanying data set are based on a real case. Only the bank's name has been changed, according to Example 11.3 of Albright, Winston and Zappe (1999).] The charge was that its female employees received substantially smaller salaries than its male employees. The bank's employee database (based on 1995 data) is listed in Albright, Winston and Zappe (1999). For each of its 208 employees the data set includes the following variables:

- EduLev: education level, a categorical variable with categories 1 (finished high school), 2 (finished some college courses), 3 (obtained a bachelor's degree), 4 (took some graduate courses), 5 (obtained a graduate degree).
- JobGrade: a categorical variable indicating the current job level, the possible levels being 1–6 (6 highest).
- YrHired: year that an employee was hired.
- YrBorn: year that an employee was born.
- Gender: a categorical variable with values "Female" and "Male."
- YrsPrior: number of years of work experience at another bank prior to working at the Fifth National Bank.
- PCJob: a dummy variable with value 1 if the empolyee's current job is computer related and value 0 otherwise.
- Salary: current (1995) annual salary in thousands of dollars.

A naive comparison of the average salaries of males and females will not work, since there are many confounding factors that affect salary. Since our main interest is to provide, after adjusting contributions from confounding factors, a good estimate for the average salary difference between male and female employees, it is very reasonable to build a large statistical model to reduce possible modeling biases. In building such a model the estimability of parameters is a factor in choosing the number of parameters, which depends on the sample size.

Two models arise naturally: the linear model

$$\text{Salary} = \beta_0 + \beta_1 \text{Female} + \beta_2 \text{PCJob} + \sum_{i=1}^{4} \beta_{2+i} \text{Edu}_i$$
$$(4.1)$$
$$+ \sum_{i=1}^{5} \beta_{6+i} \text{JobGrd}_i + \beta_{12} \text{YrsExp} + \beta_{13} \text{Age} + \varepsilon$$

and the semiparametric model

$$\text{Salary} = \beta_0 + \beta_1 \text{Female} + \beta_2 \text{PCJob} + \sum_{i=1}^{4} \beta_{2+i} \text{Edu}_i$$
$$(4.2)$$
$$+ \sum_{i=1}^{5} \beta_{6+i} \text{JobGrd}_i + f_1(\text{YrsExp}) + f_2(\text{Age}) + \varepsilon,$$

where the variable YrsExp is the years of working experience, computed from the variables YrHired and YrsPrior, and $f_1$ and $f_2$ are two continuous functions to be parameterized. Figure 1 shows the distributions of the years of working experience and age. To take into account the estimability, the number of parameters used in modeling $f_1$ and $f_2$ depends on the sample size. This falls within the framework of our study. In our analysis we employ quadratic spline models:

$$f_i(x) = \alpha_{i1} x + \alpha_{i2} x^2 + \alpha_{i3}(x - x_{i1})_+^2 + \cdots + \alpha_{i7}(x - x_{i5})_+^2, \qquad i = 1, 2,$$
(4.3)
where $x_{i1}, \ldots, x_{i5}$ are, respectively, the $2/7, 3/7, \ldots, 6/7$ sample quantiles of the variables YrsExp ($i = 1$) and Age ($i = 2$). In other words, the knots for YrsExp are 6, 8, 10, 12 and 15, and for Age are 32, 35, 37, 42 and 46. The total number of parameters for the classical linear model (4.1) is 14, and for the semiparametric model (4.2) is 26. Clearly, model (4.2) incurs smaller modeling biases than model (4.1). Since the number of parameters in both models is large, we apply the penalized likelihood method to select significant variables. Similarly to Fan and Li (2001), a modified GCV has been applied to choose the regularization parameter $\lambda$. Find $\lambda$ to minimize

$$\text{GCV}(\lambda) = \frac{1}{n} \frac{\|y - X\hat{\beta}(\lambda)\|^2}{\{1 - \gamma e(\lambda)/n\}^2},$$

where $\hat{\beta}(\lambda)$ is the penalized least-squares estimate for a given $\lambda$ and $e(\lambda)$ is the effective number of parameters defined in Section 4.2 of Fan and Li (2001). The values $\gamma = 1$ and $\gamma = 2.5$ are applied in model (4.1) and model (4.2), respectively.

There are a few cases that do not appear typical. We deleted the samples with age over 60 or working experience over 30. They correspond mainly to the company executives who earned handsome salaries. [Two of them are females who were over age 60, employed at ages 54 and 58 with no prior experience and at the lowest grade, and with just a high school education. While these two female employees had relatively low salaries, they should not be regarded as being discriminated against. Further, they have high leverages for model (4.1), particularly in the direction of age. Deleting these two cases is reasonable, either from a statistical or a practical point of view.] As a result of this deletion, a sample of size 199 remains for our analysis.

We first applied the ordinary least-squares fit. The estimated coefficients and their associated standard errors are summarized in Table 1. To apply the penalized likelihood method (1.1), we need to take care of the scale problem for each covariate. We normalized each covariate variable by the estimated standard error from the ordinary least squares and then estimated the coefficients and transformed the coefficients back to the original scale. This is equivalent to applying the penalized parameter $\lambda_j = \lambda \operatorname{SE}(\hat{\beta}_j)$ for covariate $X_j$, where $\operatorname{SE}(\hat{\beta}_j)$ is the standard error of $\hat{\beta}_j$ for the ordinary least squares estimate. The SCAD penalty function is used throughout our numerical implementation.

The penalized least-squares estimates are also presented in Table 1. For the semiparametric model (4.2) the regression components for the variables YrsExp and Age are shown in Figure 2. The residual plots against the variables YrsExp and Age are shown in Figure 3. They do not exhibit any systematic patterns.

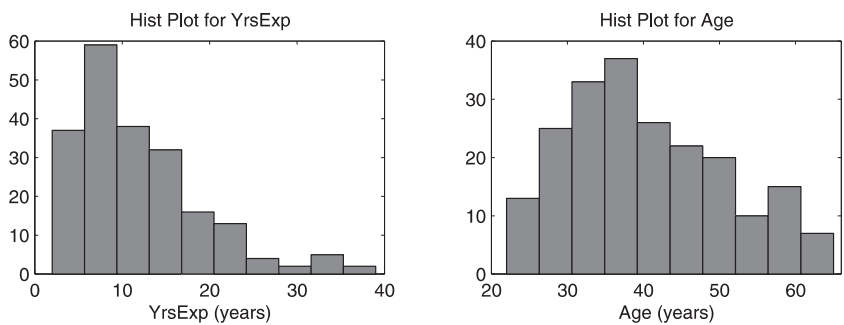

Fig. 1. *Distributions of years of working experience and age.*

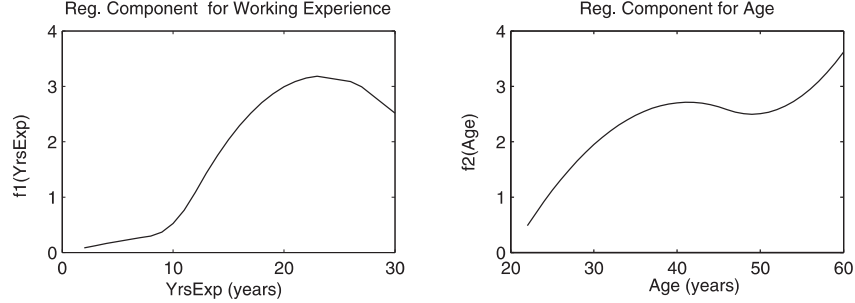

Fig. 2. *Regression components $f_1$ and $f_2$ for the semiparametric model (4.2).*

The three models all have high multiple $R^2$—over 80% of the salary variation can be explained by the variables that we use. None of the results shows statistical evidence of discrimination. The coefficients in front of the indicator of Female are negative, but statistically insignificant.

We now apply the likelihood ratio test to examine whether there is any discrimination against female employees. This leads to the following null hypothesis:

$$H_0 : \beta_1 = 0.$$

Even in the presence of a large number of parameters, according to Theorem 4, the likelihood ratio theory continues to apply. Table 2 summarizes the test results under both models (4.1) and (4.2), using both penalized and unpenalized versions of the likelihood ratio test. The results are consistent with the regression analyses depicted in Table 1.

Table 1
*Estimates and standard errors for Fifth National Bank data*

| Method    | Least squares |         | SCAD PLS |         | SCAD PLS |         |
|-----------|---------------|---------|----------|---------|----------|---------|
| Intercept | 54.238        | (2.067) | 55.835   | (1.527) | 52.470   | (2.890) |
| Female    | −0.556        | (0.637) | −0.624   | (0.639) | −0.933   | (0.708) |
| PcJob     | 3.982         | (0.908) | 4.151    | (0.909) | 2.851    | (0.640) |
| Ed1       | −1.739        | (1.049) | 0        | (—)     | 0        | (—)     |
| Ed2       | −2.866        | (0.999) | −1.074   | (0.522) | −0.542   | (0.265) |
| Ed3       | −2.145        | (0.753) | −0.914   | (0.421) | 0        | (—)     |
| Ed4       | −1.484        | (1.369) | 0        | (—)     | 0        | (—)     |
| Job1      | −22.954       | (1.734) | −24.643  | (1.535) | −22.841  | (1.332) |
| Job2      | −21.388       | (1.686) | −22.818  | (1.546) | −20.591  | (1.370) |
| Job3      | −17.642       | (1.634) | −18.803  | (1.562) | −16.719  | (1.391) |
| Job4      | −13.046       | (1.578) | −13.859  | (1.529) | −11.807  | (1.359) |
| Job5      | −7.462        | (1.551) | −7.770   | (1.539) | −5.235   | (1.150) |
| YrsExp    | 0.215         | (0.065) | 0.193    | (0.046) | (—)      | (—)     |
| Age       | 0.030         | (0.039) | 0        | (—)     | (—)      | (—)     |
| $R^2$     | 0.8221        |         | 0.8176   |         | 0.8123   |         |

TABLE 2
*SCAD penalized likelihood ratio test*

|  | Likelihood ratio test | | Penalized likelihood ratio test | |
| --- | --- | --- | --- | --- |
|  | $\chi^2$-statistic | $P$-value | $\chi^2$-statistic | $P$-value |
| Model (1.1) | 0.7607 | 0.3831 | 1.0131 | 0.3142 |
| Model (1.2) | 1.8329 | 0.1758 | 1.4311 | 0.2316 |

One can apply the logarithmic transformation to the salary variable and analyze the transformed data. This would make the transformed data more normally distributed. We opted for the original scale for the sake of interpretability. Furthermore, the conclusion does not change much.

While there is no significant statistical evidence for discrimination based on the above analyses, the arguments can still go on. For example, as intuitively expected, the job grade is a very important variable that determines the salary. For this data set, it explains 77.29% of the salary variation. Now the question arises naturally whether it was harder for females employees to be promoted, after adjusting for variables such as working experience, age and education level. We do not pursue this issue further.

4.2. *A simulation study.* In this section we use a simulation study to augment our theoretical results. To present a situation in which the number of parameters depends on $n$, to show the applicability of our results is wider than what we have presented and to create dependence between covariates, we consider the following autoregressive model:

(4.4) $\quad X_i = \beta_1 X_{i-1} + \beta_2 X_{i-2} + \cdots + \beta_p X_{i-p_n} + \varepsilon, \qquad i = 1, 2, \ldots, n,$

where $\beta = (11/4, -23/6, 37/12, -13/9, 1/3, 0, \ldots, 0)^T$ and $\varepsilon$ is white noise with variance $\sigma^2$. The number of parameters depends naturally on $n$, as time series analysts naturally wish to explore the order of fit to reduce possible modeling biases. This time series model is stationary since its associated polynomial

$$\Phi(B) = \left(1 - \frac{3B}{4}\right)\left(1 - B + \frac{2B^2}{3}\right)^2$$

has no zero inside the unit circle.

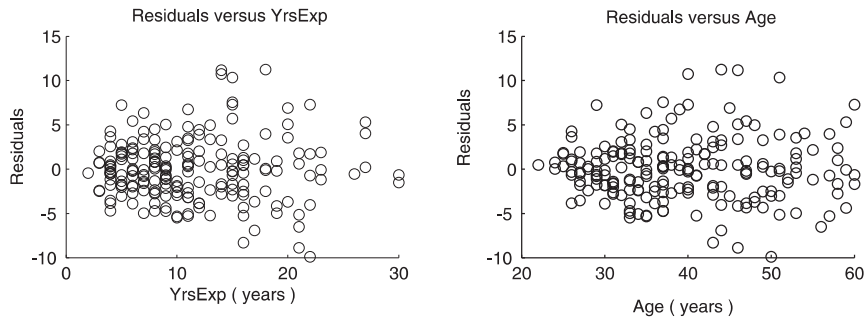

FIG. 3. *Residuals after fitting the semiparametric model (4.2).*

In our simulation experiments 400 samples of sizes 100, 200, 400 and 800 with $p_n = [4n^{1/4}] - 5$ are drawn from model (4.4). The penalized least-squares method with the SCAD penalty is employed. The algorithm of Fan and Li (2001) is used. The medians of relative model errors (MRMEs) $(\hat{\beta} - \beta)^T E X^T X (\hat{\beta} - \beta)$, among the least-squares (LS) estimator, the penalized least-squares (PLS) estimator and the oracle estimator, measure the effectiveness of the methods. The results are summarized in Table 3. As expected, the oracle estimator performs the best and the PLS performs comparably with the oracle estimator for all sample sizes. The LS estimator performs the worst and its relative performance deteriorates as $n$ increases. The average number of zero coefficients is also reported in Table 3, in which the column labeled "Correct" presents the average number restricted only to the true zero coefficients, and the column labeled "Incorrect" depicts the average number of coefficients erroneously set to 0. For example, for $n = 400$, among seven nonzero coefficients, on average 5.78 coefficients, or 83%, were correctly estimated as 0, and among five nonzero coefficients, on average 0.22 coefficient was incorrectly estimated as 0. The medians of the estimated coefficients over 400 simulated data sets are presented in Table 4. The biases are quite small, except for estimated $\beta_5$, which is difficult to estimate. The variances of the estimated coefficients across 400 simulations are presented in Table 5.

To test the accuracy of the standard error formula, the standard deviations of the estimated coefficients are computed among 400 simulations. These can be regarded as the true standard errors (columns labeled SD in

TABLE 3
*Simulation results for the time series model*

| $n$ | $p_n$ | MRME (%) | | | Average number of zero coeficients | |
|---|---|---|---|---|---|---|
| | | Oracle/LS | PLS/LS | Oracle/PLS | Correct | Incorrect |
| 100 | 7 | 75.33 | 89.21 | 80.17 | 1.34 [67%] | 0.49 |
| 200 | 10 | 50.61 | 69.64 | 73.27 | 3.91 [78%] | 0.39 |
| 400 | 12 | 40.03 | 59.57 | 73.06 | 5.78 [83%] | 0.22 |
| 800 | 16 | 31.75 | 49.05 | 70.08 | 9.49 [86%] | 0.10 |

TABLE 4
*Median of estimators for coefficients of time series model*

| $n$ | $p_n$ | $\hat{\beta}_1$ | $\hat{\beta}_2$ | $\hat{\beta}_3$ | $\hat{\beta}_4$ | $\hat{\beta}_5$ |
|---|---|---|---|---|---|---|
| 100  | 7  | 2.678 | −3.616 | 2.739 | −1.096 | 0 |
| 200  | 10 | 2.711 | −3.696 | 2.856 | −1.240 | 0.242 |
| 400  | 12 | 2.729 | −3.769 | 2.959 | −1.333 | 0.293 |
| 800  | 16 | 2.737 | −3.792 | 3.023 | −1.383 | 0.306 |
| True | —  | 2.750 | −3.833 | 3.083 | −1.444 | 0.333 |

Table 5) and compared with 400 estimated standard errors. The 400 estimated standard errors are summarized by their median (columns $SD_m$) and interquartile range divided by 1.349 (columns $SD_{mad}$), which is a robust estimate of the standard deviation. The results are presented in Table 5. The accuracy gets better when $n$ increases. Further, the accuracy is better for the first two coefficients than for the last two coefficients, which have lower signal-to-noise ratios.

We now apply the penalized likelihood ratio test to the following null hypothesis:

$$H_0 : \beta_6 = \beta_7 = 0.$$

The likelihood ratio statistic is computed for each simulation. The distribution of these statistics among 400 simulations can be regarded as the true null distribution and can be compared with the asymptotic distribution. Figure 4 depicts the estimated densities of the likelihood ratio statistics among 400 simulations for $n = 100$ and 400.

**5. Proofs of theorems.** In this section, we give rigorous proofs of Theorems 1–4.

PROOF OF THEOREM 1. Let $\alpha_n = \sqrt{p_n}(n^{-1/2} + a_n)$ and set $\|\mathbf{u}\| = C$, where $C$ is a large enough constant. Our aim is to show that for any given

TABLE 5
*Standard deviations (multiplied by* 1000*) of estimators for time series model*

| | $\hat{\beta}_1$ | | $\hat{\beta}_2$ | | $\hat{\beta}_3$ | | $\hat{\beta}_4$ | | $\hat{\beta}_5$ | |
|---|---|---|---|---|---|---|---|---|---|---|
| | | $SD_m$ | | $SD_m$ | | $SD_m$ | | $SD_m$ | | $SD_m$ |
| $n$ | SD | ($SD_{mad}$) | SD | ($SD_{mad}$) | SD | ($SD_{mad}$) | SD | ($SD_{mad}$) | SD | ($SD_{mad}$) |
| 100 | 120 | 91 (5.1) | 337 | 230 (29.8) | 525 | 285 (66.6) | 451 | 177 (87.2) | 249 | 79 (66.7) |
| 200 | 76  | 66 (2.8) | 221 | 174 (15.2) | 340 | 231 (58)   | 348 | 170 (87.2) | 243 | 64 (49.5) |
| 400 | 50  | 47 (1.2) | 149 | 126 (4.5)  | 222 | 169 (8.8)  | 204 | 125 (9.0)  | 129 | 47 (3.90) |
| 800 | 35  | 34 (0.7) | 99  | 90 (3.1)   | 145 | 121 (8.5)  | 132 | 90 (14.1)  | 63  | 34 (12.5) |

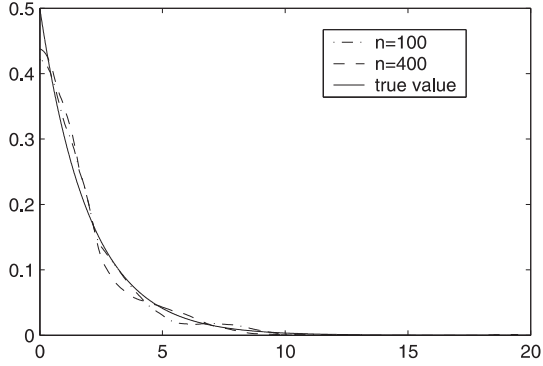

FIG. 4. *Estimated densities of the likelihood ratio statistics for $n = 100$ (dot-dash) and $n = 400$ (long-dash) along with the density of the $\chi_2^2$ distribution (solid).*

$\varepsilon$ there is a large constant $C$ such that, for large $n$ we have

(5.1) $$P\left\{\sup_{\|\mathbf{u}\|=C} Q_n(\beta_{n0} + \alpha_n \mathbf{u}) < Q_n(\beta_{n0})\right\} \geq 1 - \varepsilon.$$

This implies that with probability tending to 1 there is a local maximum $\hat{\beta}_n$ in the ball $\{\beta_{n0} + \alpha_n \mathbf{u} : \|\mathbf{u}\| \leq C\}$ such that $\|\hat{\beta}_n - \beta_{n0}\| = O_p(\alpha_n)$.

Using $p_{\lambda_n}(0) = 0$, we have

$$D_n(\mathbf{u}) = Q_n(\beta_{n0} + \alpha_n \mathbf{u}) - Q_n(\beta_{n0})$$
$$\leq L_n(\beta_{n0} + \alpha_n \mathbf{u}) - L_n(\beta_{n0})$$
$$- n \sum_{j=1}^{s_n} \{p_{\lambda_n}(|\beta_{n0j} + \alpha_n u_j|) - p_{\lambda_n}(|\beta_{n0j}|)\}$$
$$\triangleq (I) + (II).$$

Then by Taylor's expansion we obtain

$$(I) = \alpha_n \nabla^T L_n(\beta_{n0})\mathbf{u} + \tfrac{1}{2}\mathbf{u}^T \nabla^2 L_n(\beta_{n0})\mathbf{u}\alpha_n^2 + \tfrac{1}{6}\nabla^T\{\mathbf{u}^T \nabla^2 L_n(\beta_n^*)\mathbf{u}\}\mathbf{u}\alpha_n^3$$
$$\triangleq I_1 + I_2 + I_3,$$

where the vector $\beta_n^*$ lies between $\beta_{n0}$ and $\beta_{n0} + \alpha_n \mathbf{u}$, and

$$(II) = -\sum_{j=1}^{s_n}[n\alpha_n p'_{\lambda_n}(|\beta_{n0j}|)\operatorname{sgn}(\beta_{n0j})u_j + n\alpha_n^2 p''_n(\beta_{n0j})u_j^2\{1 + o(1)\}]$$
$$\triangleq I_4 + I_5.$$

By condition (F),

(5.2) $$|I_1| = |\alpha_n \nabla^T L_n(\beta_{n0})\mathbf{u}| \leq \alpha_n \|\nabla^T L_n(\beta_{n0})\| \|\mathbf{u}\|$$
$$= O_p(\alpha_n \sqrt{np_n})\|\mathbf{u}\| = O_p(\alpha_n^2 n)\|\mathbf{u}\|.$$

Next we consider $I_2$. An application of Lemma 8 in the Appendix yields that

$$I_2 = \frac{1}{2}\mathbf{u}^T \left[\frac{1}{n}\{\nabla^2 L_n(\beta_{n0}) - E\nabla^2 L_n(\beta_{n0})\}\right]\mathbf{u} \cdot n\alpha_n^2$$

(5.3)
$$- \frac{1}{2}\mathbf{u}^T I_n(\beta_{n0})\mathbf{u} \cdot n\alpha_n^2$$

$$= -\frac{n\alpha_n^2}{2}\mathbf{u}^T I_n(\beta_{n0})\mathbf{u} + o_p(1) \cdot n\alpha_n^2 \|\mathbf{u}\|^2.$$

By the Cauchy–Schwarz inequality and condition (G), we have

$$|I_3| = \left|\frac{1}{6}\sum_{i,j,k=1}^{p_n} \frac{\partial L_n(\beta_n^*)}{\partial \beta_{ni}\, \partial \beta_{nj}\, \partial \beta_{nk}} u_i u_j u_k \alpha_n^3\right|$$

$$\leq \frac{1}{6}\sum_{l=1}^{n}\left\{\sum_{i,j,k=1}^{p_n} M_{nijk}^2(V_{nl})\right\}^{1/2} \|\mathbf{u}\|^3 \alpha_n^3.$$

Since $p_n^4/n \to 0$ and $p_n^2 a_n \to 0$ as $n \to \infty$, we have

$$\tfrac{1}{6}\sum_{l=1}^{n}\left\{\sum_{i,j,k=1}^{p_n} M_{nijk}^2(V_{nl})\right\}^{1/2} \|\mathbf{u}\|^3 \alpha_n^3$$

$$= O_p(p_n^{3/2}\alpha_n) n\alpha_n^2 \|\mathbf{u}\|^2 = o_p(n\alpha_n^2)\|\mathbf{u}\|^2.$$

Thus,

(5.4) $$I_3 = o_p(n\alpha_n^2)\|\mathbf{u}\|^2.$$

The terms $I_4$ and $I_5$ can be dealt with as follows. First,

(5.5) $$|I_4| \leq \sum_{j=1}^{s_n} |n\alpha_n p'_{\lambda_n}(|\beta_{n0j}|)\operatorname{sgn}(\beta_{n0j}) u_j| \leq \sqrt{s_n} \cdot n\alpha_n a_n \|\mathbf{u}\| \leq n\alpha_n^2 \|\mathbf{u}\|$$

and

(5.6) $$I_5 = \sum_{j=1}^{s_n} n\alpha_n^2 p''_{\lambda_n}(|\beta_{n0j}|) u_j^2 \{1+o(1)\} \leq 2 \cdot \max_{1\leq j \leq s_n} p''_{\lambda_n}(|\beta_{n0j}|) \cdot n\alpha_n^2 \|\mathbf{u}\|^2.$$

By (5.2)–(5.6) and condition (C), and allowing $\|\mathbf{u}\|$ to be large enough, all terms $I_1, I_3, I_4$ and $I_5$ are dominated by $I_2$, which is negative. This proves (5.1). $\square$

To prove Theorem 2, we first show that the nonconcave penalized estimator possesses the sparsity property $\hat{\beta}_{n2} = 0$ by the following lemma.

LEMMA 5. *Assume that conditions* (A) *and* (E)–(H) *are satisfied. If* $\lambda_n \to 0$, $\sqrt{n/p_n}\lambda_n \to \infty$ *and* $p_n^5/n \to 0$ *as* $n \to \infty$, *then with probability tending to* 1, *for any given* $\beta_{n1}$ *satisfying* $\|\beta_{n1} - \beta_{n01}\| = O_p(\sqrt{p_n/n})$ *and any constant* $C$,

$$Q\{(\beta_{n1}^T, 0)^T\} = \max_{\|\beta_{n2}\| \leq C(p_n/n)^{1/2}} Q\{(\beta_{n1}^T, \beta_{n2}^T)^T\}.$$

PROOF. Let $\varepsilon_n = C\sqrt{p_n/n}$. It is sufficient to show that with probability tending to 1 as $n \to \infty$, for any $\beta_{n1}$ satisfying $\beta_{n1} - \beta_{n01} = O_p(\sqrt{p_n/n})$ we have, for $j = s_n + 1, \ldots, p_n$,

(5.7) $$\frac{\partial Q_n(\beta_n)}{\partial \beta_{nj}} < 0 \quad \text{for } 0 < \beta_{nj} < \varepsilon_n,$$

(5.8) $$\frac{\partial Q_n(\beta_n)}{\partial \beta_{nj}} > 0 \quad \text{for } -\varepsilon_n < \beta_{nj} < 0.$$

By Taylor expansion,

$$\frac{\partial Q_n(\beta_n)}{\partial \beta_{nj}} = \frac{\partial L_n(\beta_n)}{\partial \beta_{nj}} - np'_{\lambda_n}(|\beta_{nj}|)\operatorname{sgn}(\beta_{nj})$$

$$= \frac{\partial L_n(\beta_{n0})}{\partial \beta_{nj}} + \sum_{l=1}^{p_n} \frac{\partial^2 L_n(\beta_{n0})}{\partial \beta_{nj} \partial \beta_{nl}}(\beta_{nl} - \beta_{n0l})$$

$$+ \sum_{l,k=1}^{p_n} \frac{\partial^3 L_n(\beta_n^*)}{\partial \beta_{nj} \partial \beta_{nl} \partial \beta_{nk}}(\beta_{nl} - \beta_{n0l})(\beta_{nk} - \beta_{n0k})$$

$$- np'_{\lambda_n}(|\beta_{nj}|)\operatorname{sgn}(\beta_{nj})$$

$$\widehat{=} I_1 + I_2 + I_3 + I_4,$$

where $\beta_n^*$ lies between $\beta_n$ and $\beta_{n0}$. Next, we consider $I_1$, $I_2$ and $I_3$.

By a standard argument, we have

(5.9) $$I_1 = O_p(\sqrt{n}) = O_p(\sqrt{np_n}).$$

The term $I_2$ can be written as

$$I_2 = \sum_{l=1}^{p_n} \left\{ \frac{\partial^2 L_n(\beta_{n0})}{\partial \beta_{nj} \partial \beta_{nl}} - E\frac{\partial^2 L_n(\beta_{n0})}{\partial \beta_{nj} \partial \beta_{nl}} \right\}(\beta_{nl} - \beta_{n0l})$$

$$+ \sum_{l=1}^{p_n} E\frac{\partial^2 L_n(\beta_{n0})}{\partial \beta_{nj} \partial \beta_{nl}}(\beta_{nl} - \beta_{n0l})$$

$$\widehat{=} K_1 + K_2.$$

Using the Cauchy–Schwarz inequality and $\|\beta_n - \beta_{n0}\| = O_p(\sqrt{p_n/n})$, we have

$$|K_2| = \left| n \sum_{l=1}^{p_n} I_n(\beta_{n0})(j,l)(\beta_{nl} - \beta_{n0l}) \right|$$

$$\leq nO_p\left(\sqrt{\frac{p_n}{n}}\right) \left\{ \sum_{l=1}^{p_n} I_n^2(\beta_{n0})(j,l) \right\}^{1/2}.$$

As the eigenvalues of the information matrix are bounded according to condition (F), we have

$$\sum_{l=1}^{p_n} I_n^2(\beta_{n0})(j,l) = O(1).$$

This entails that

(5.10) $$K_2 = O_p(\sqrt{np_n}).$$

As for the term $K_1$, by the Cauchy–Schwarz inequality we have

$$|K_1| \leq \|\beta_n - \beta_{n0}\| \left[ \sum_{l=1}^{p_n} \left\{ \frac{\partial^2 L_n(\beta_{n0})}{\partial \beta_{nj} \partial \beta_{nl}} - E \frac{\partial^2 L_n(\beta_{n0})}{\partial \beta_{nj} \partial \beta_{nl}} \right\}^2 \right]^{1/2}.$$

Then from condition (F), it is easy to show that

$$\left[ \sum_{l=1}^{p_n} \left\{ \frac{\partial^2 L_n(\beta_{n0})}{\partial \beta_{nj} \partial \beta_{nl}} - E \frac{\partial^2 L_n(\beta_{n0})}{\partial \beta_{nj} \partial \beta_{nl}} \right\}^2 \right]^{1/2} = O_p(\sqrt{np_n}).$$

By $\|\beta_n - \beta_{n0}\| = O_p(\sqrt{p_n/n})$ it follows that $K_1 = O_p(\sqrt{np_n})$. This, together with (5.10), yields

(5.11) $$I_2 = O_p(\sqrt{np_n}).$$

Next we consider $I_3$. We can write $I_3$ as follows:

$$I_3 = \sum_{l,k=1}^{p_n} \left\{ \frac{\partial^3 L_n(\beta_n^*)}{\partial \beta_{nj} \beta_{nl} \beta_{nk}} - E \frac{\partial^3 L_n(\beta_n^*)}{\partial \beta_{nj} \beta_{nl} \beta_{nk}} \right\} (\beta_{nj} - \beta_{n0j})(\beta_{nk} - \beta_{n0k})$$

$$+ \sum_{l,k=1}^{p_n} E \frac{\partial^3 L_n(\beta_n^*)}{\partial \beta_{nj} \beta_{nl} \beta_{nk}} (\beta_{nj} - \beta_{n0j})(\beta_{nk} - \beta_{n0k})$$

$$\hat{=} K_3 + K_4.$$

By condition (G),

(5.12) $$|K_4| \leq C_5^{1/2} \cdot np_n \cdot \|\beta_n - \beta_{n0}\|^2 = O_p(p_n^2) = o_p(\sqrt{np_n}).$$

However, by the Cauchy–Schwarz inequality,

$$K_3^2 \leq \sum_{l,k=1}^{p_n} \left\{ \frac{\partial^3 L_n(\beta_n^*)}{\partial \beta_{nj} \partial \beta_{nk} \partial \beta_{nl}} - E\frac{\partial^3 L_n(\beta_n^*)}{\partial \beta_{nj} \partial \beta_{nk} \partial \beta_{nl}} \right\}^2 \|\beta_n - \beta_{n0}\|^4.$$

Under conditions (G) and (H), we have

(5.13) $$K_3 = O_p\left\{ \left(np_n^2 \frac{p_n^2}{n^2}\right)^{1/2} \right\} = o_p(\sqrt{np_n}).$$

From (5.9) and (5.11)–(5.13) we have

$$I_1 + I_2 + I_3 = O_p(\sqrt{np_n}).$$

Because $\sqrt{p_n/n}/\lambda_n \to 0$ and $\liminf_{n\to\infty} \inf_{\theta \to 0^+} p'_{\lambda_n}(\theta)/\lambda_n > 0$, from

$$\frac{\partial Q_n(\beta_n)}{\partial \beta_{nj}} = n\lambda_n \left\{ -\frac{p'_{\lambda_n}(|\beta_{nj}|)}{\lambda_n} \operatorname{sgn}(\beta_{nj}) + O_p\left(\sqrt{\frac{p_n}{n}}/\lambda_n\right) \right\}$$

it is easy to see that the sign of $\beta_{nj}$ completely determines the sign of $\partial Q_n(\beta_n)/\partial \beta_{nj}$. Hence, (5.7) and (5.8) follow. □

PROOF OF THEOREM 2. As shown in Theorem 1, there is a root-$(n/p_n)$-consistent local maximizer $\hat{\beta}_n$ of $Q_n(\beta_n)$. By Lemma 5, part (i) holds that $\hat{\beta}_n$ has the form $(\hat{\beta}_{n1}, 0)^T$. We need only prove part (ii), the asymptotic normality of the penalized nonconcave likelihood estimator $\hat{\beta}_{n1}$.

If we can show that

$$\{I_n(\beta_{n01}) + \Sigma_{\lambda_n}\}(\hat{\beta}_{n1} - \beta_{n01}) + \mathbf{b}_n = \frac{1}{n}\nabla L_n(\beta_{n01}) + o_p(n^{-1/2}),$$

then

$$\sqrt{n}A_n I_n^{-1/2}(\beta_{n01})\{I_n(\beta_{n01}) + \Sigma_{\lambda_n}\}[\hat{\beta}_{n1} - \beta_{n01} + \{I_n(\beta_{n01}) + \Sigma_{\lambda_n}\}^{-1}\mathbf{b}_n]$$
$$= \frac{1}{\sqrt{n}}A_n I_n^{-1/2}(\beta_{n01})\nabla L_n(\beta_{n01}) + o_p\{A_n I_n^{-1/2}(\beta_{n01})\}.$$

By the conditions of Theorem 2, we have the last term of $o_p(1)$. Let

$$Y_{ni} = \frac{1}{\sqrt{n}}A_n I_n^{-1/2}(\beta_{n01})\nabla L_{ni}(\beta_{n01}), \qquad i = 1, 2, \ldots, n.$$

It follows that, for any $\varepsilon$,

$$\sum_{i=1}^n E\|Y_{ni}\|^2 \mathbf{1}\{\|Y_{ni}\| > \varepsilon\} = nE\|Y_{n1}\|^2 \mathbf{1}\{\|Y_{n1}\| > \varepsilon\}$$
$$\leq n\{E\|Y_{n1}\|^4\}^{1/2}\{P(\|Y_{n1}\| > \varepsilon)\}^{1/2}.$$

By condition (F) and $A_n A_n^T \to G$, we obtain

$$P(\|Y_{n1}\| > \varepsilon) \leq \frac{E\|A_n I_n^{-1/2}(\beta_{n01}) \nabla L_{n1}(\beta_{n01})\|^2}{n\varepsilon} = O(n^{-1})$$

and

$$E\|Y_{n1}\|^4 = \frac{1}{n^2} E\|A_n I_n^{-1/2}(\beta_{n01}) \nabla L_{n1}(\beta_{n01})\|^4$$

$$\leq \frac{1}{n^2} \lambda_{\max}(A_n A_n^T) \lambda_{\max}\{I_n(\beta_{n01})\} E\|\nabla^T L_{n1}(\beta_{n01}) \nabla L_{n1}(\beta_{n01})\|^2$$

$$= O\left(\frac{p_n^2}{n^2}\right).$$

Thus, we have

$$\sum_{i=1}^n E\|Y_{ni}\|^2 \mathbf{1}\{\|Y_{ni}\| > \varepsilon\} = O\left(n \frac{p_n}{n} \frac{1}{\sqrt{n}}\right) = o(1).$$

On the other hand, as $A_n A_n^T \to G$, we have

$$\sum_{i=1}^n \text{cov}(Y_{ni}) = n \text{cov}(Y_{n1}) = \text{cov}\{A_n I_n^{-1/2}(\beta_{n01}) \nabla L_{n1}(\beta_{n01})\} \to G.$$

Thus, $Y_{ni}$ satisfies the conditions of the Lindeberg–Feller central limit theorem [see van der Vaart (1998)]. This also means that $1/\sqrt{n} A_n I_n(\beta_{n01})^{-1/2} \nabla L_n(\beta_{n01})$ has an asymptotic multivariate normal distribution.

With a slight abuse of notation, let $Q_n(\beta_{n1}) = Q_n(\beta_{n1}, 0)$. As $\hat{\beta}_{n1}$ must satisfy the penalized likelihood equation $\nabla Q_n(\hat{\beta}_{n1}) = 0$, using the Taylor expansion on $\nabla Q_n(\hat{\beta}_{n1})$ at point $\beta_{n01}$, we have

$$\frac{1}{n}[\{\nabla^2 L_n(\beta_{n01}) - \nabla^2 P_{\lambda_n}(\beta_{n1}^{**})\}(\hat{\beta}_{n1} - \beta_{n01}) - \nabla P_{\lambda_n}(\beta_{n01})]$$

$$= -\frac{1}{n}\left[\nabla L_n(\beta_{n01}) + \frac{1}{2}(\hat{\beta}_{n1} - \beta_{n01})^T \nabla^2\{\nabla L_n(\beta_{n1}^*)\}(\hat{\beta}_{n1} - \beta_{n01})\right],$$

where $\beta_{n1}^*$ and $\beta_{n1}^{**}$ lie between $\hat{\beta}_{n1}$ and $\beta_{n01}$. Now we define

$$\mathcal{L} \triangleq \nabla^2 L_n(\beta_{n01}) - \nabla^2 P_{\lambda_n}(\beta_{n1}^{**})$$

and

$$\mathcal{C} \triangleq \tfrac{1}{2}(\hat{\beta}_{n1} - \beta_{n01})^T \nabla^2\{\nabla L_n(\beta_{n1}^*)\}(\hat{\beta}_{n1} - \beta_{n01}).$$

Under conditions (G) and (H) and by the Cauchy–Schwarz inequality, we have

(5.14)
$$\left\|\frac{1}{n}\mathcal{C}\right\|^2 \leq \frac{1}{n^2} \sum_{i=1}^n n^2 \|\hat{\beta}_{n1} - \beta_{n01}\|^4 \sum_{j,k,l=1}^{s_n} M_{njkl}^2(V_{ni})$$

$$= O_p\left(\frac{p_n^2}{n^2}\right) O_p(p_n^3) = o_p\left(\frac{1}{n}\right).$$

At the same time, by Lemma 8 in the Appendix and because of condition (H), it is easy to show that

$$\lambda_i\left\{\frac{1}{n}\mathcal{L} + I_n(\beta_{n01}) + \Sigma_{\lambda_n}\right\} = o_p\left(\frac{1}{\sqrt{p_n}}\right), \qquad i = 1, 2, \ldots, s_n,$$

where $\lambda_i(A)$ is the $i$th eigenvalue of a symmetric matrix $A$. As $\hat{\beta}_{n1} - \beta_{n01} = O_p(\sqrt{p_n/n})$,

(5.15) $$\left\{\frac{1}{n}\mathcal{L} + I_n(\beta_{n01}) + \Sigma_{\lambda_n}\right\}(\hat{\beta}_{n1} - \beta_{n01}) = o_p\left(\frac{1}{\sqrt{n}}\right).$$

Finally, from (5.14) and (5.15) we have

(5.16) $$\{I_n(\beta_{n01}) + \Sigma_{\lambda_n}\}(\hat{\beta}_{n1} - \beta_{n01}) + \mathbf{b}_n = \frac{1}{n}\nabla L_n(\beta_{n01}) + o_p\left(\frac{1}{\sqrt{n}}\right).$$

Following (5.16), Theorem 2 follows. □

PROOF OF THEOREM 3. Let $\mathcal{A}_n = -n^{-1}\nabla^2 L_n(\hat{\beta}_{n1}) + \Sigma_{\lambda_n}(\hat{\beta}_{n1})$, $\mathcal{B}_n = \widehat{\text{cov}}\{\nabla L_n(\hat{\beta}_{n1})\}$, $\mathcal{A} = I_n(\beta_{n01}) + \Sigma_{\lambda_n}$ and $\mathcal{B} = I_n(\beta_{n01})$. Then we have

$$\hat{\Sigma}_n - \Sigma_n = \mathcal{A}_n^{-1}(\mathcal{B}_n - \mathcal{B})\mathcal{A}_n^{-1} + (\mathcal{A}_n^{-1} - \mathcal{A}^{-1})\mathcal{B}\mathcal{A}_n^{-1} + \mathcal{A}^{-1}\mathcal{B}(\mathcal{A}_n^{-1} - \mathcal{A}^{-1})$$
$$= I_1 + I_2 + I_3$$

and

$$\mathcal{A}_n^{-1} - \mathcal{A}^{-1} = \mathcal{A}_n^{-1}(\mathcal{A} - \mathcal{A}_n)\mathcal{A}^{-1}.$$

Let $\lambda_i(\mathcal{A})$ be the $i$th eigenvalue of a symmetric matrix $\mathcal{A}$. If we can show that $\lambda_i(\mathcal{A} - \mathcal{A}_n) = o_p(1)$ and $\lambda_i(\mathcal{B}_n - \mathcal{B}) = o_p(1)$, then from the fact that $|\lambda_i(\mathcal{B})|$ and $|\lambda_i(\mathcal{A})|$ are uniformly bounded away from 0 and infinite, we have

$$\lambda_i(\hat{\Sigma}_n - \Sigma_n) = o_p(1).$$

This means that $\hat{\Sigma}_n$ is a weakly consistent estimator of $\Sigma_n$.

First, let us consider $\mathcal{A} - \mathcal{A}_n$ and decompose it as follows:

$$\mathcal{A} - \mathcal{A}_n = I_n(\beta_{n01}) + \frac{1}{n}\nabla^2 L_n(\hat{\beta}_{n1}) + \Sigma_{\lambda_n}(\beta_{n01}) - \Sigma_{\lambda_n}(\hat{\beta}_{n1}) \hat{=} K_1 + K_2.$$

It is obvious that

$$\lambda_{\min}(K_1) + \lambda_{\min}(K_2) \leq \lambda_{\min}(K_1 + K_2)$$
$$\leq \lambda_{\max}(K_1 + K_2) \leq \lambda_{\max}(K_1) + \lambda_{\max}(K_2).$$

Thus, we need only consider $\lambda_i(K_1)$ and $\lambda_i(K_2)$ separately. The term $K_1$ can be expressed as

$$K_1 = I_n(\beta_{n01}) + \frac{1}{n}\nabla^2 L_n(\beta_{n01}) - \frac{1}{n}\nabla^2 L_n(\beta_{n01}) + \frac{1}{n}\nabla^2 L_n(\hat{\beta}_{n1}).$$

According to Lemma 8 in the Appendix, we have

$$\left\| I_n(\beta_{n01}) + \frac{1}{n}\nabla^2 L_n(\beta_{n01}) \right\| = o_p(1). \tag{5.17}$$

As shown in Lemma 9,

$$\left\| \frac{\nabla^2 L_n(\hat{\beta}_{n1})}{n} - \frac{\nabla^2 L_n(\beta_{n01})}{n} \right\|^2 = O_p\left(\frac{p_n^4}{n}\right) = o_p(1). \tag{5.18}$$

Thus, it follows from (5.17) and (5.18) that $\|K_1\| = o_p(1)$. This also means that we have

$$\lambda_i(K_1) = o_p(1), \qquad i = 1, 2, \ldots, s_n. \tag{5.19}$$

As $\|\hat{\beta}_{n1} - \beta_{n01}\| = O_p(\sqrt{p_n/n})$, by condition (D), $p''_{\lambda_n}(\hat{\beta}_{nj}) - p''_{\lambda_n}(\beta_{n0j}) = o_p(1)$, that is,

$$\lambda_i(K_2) = o_p(1), \qquad i = 1, 2, \ldots, s_n. \tag{5.20}$$

Hence, from (5.19) and (5.20) we have shown that

$$\lambda_i(\mathcal{A} - \mathcal{A}_n) = o_p(1), \qquad i = 1, 2, \ldots, s_n. \tag{5.21}$$

Next we consider $\lambda_i(\mathcal{B}_n - \mathcal{B})$. First we express $\mathcal{B}_n - \mathcal{B}$ as the sum of $K_3$ and $K_4$, where

$$K_3 \triangleq \left\{ \frac{1}{n} \sum_{i=1}^n \frac{\partial L_{ni}(\hat{\beta}_{n1})}{\partial \beta_j} \frac{\partial L_{ni}(\hat{\beta}_{n1})}{\partial \beta_k} \right\} - I_n(\beta_{n01})$$

and

$$K_4 \triangleq -\left\{ \frac{1}{n} \sum_{i=1}^n \frac{\partial L_{ni}(\hat{\beta}_{n1})}{\partial \beta_j} \right\} \left\{ \frac{1}{n} \sum_{i=1}^n \frac{\partial L_{ni}(\hat{\beta}_{n1})}{\partial \beta_k} \right\}.$$

Using the aforementioned argument, we need only consider $K_3$ and $K_4$ separately.

Note that

$$\frac{1}{n} \sum_{i=1}^n \frac{\partial L_{ni}(\hat{\beta}_{n1})}{\partial \beta_j} - p'_{\lambda_n}(|\hat{\beta}_{nj}|) = 0, \qquad j = 1, 2, \ldots, s_n,$$

which implies that

$$\|K_4\|^2 = \sum_{j=1}^{s_n} \sum_{k=1}^{s_n} \{p'_{\lambda_n}(|\hat{\beta}_{nj}|)\}^2 \{p'_{\lambda_n}(|\hat{\beta}_{nk}|)\}^2$$

$$= \left\{ \sum_{j=1}^{s_n} p'_{\lambda_n}(|\hat{\beta}_{nj}|)^2 \right\}^2. \tag{5.22}$$

By Taylor expansion,

(5.23) $$p'_{\lambda_n}(|\hat{\beta}_{nj}|) = p'_{\lambda_n}(|\beta_{n0j}|) + p''_{\lambda_n}(|\beta^*_{nj}|)(\hat{\beta}_{nj} - \beta_{n0j}),$$

where $\beta^*_{nj}$ lies between $\hat{\beta}_{nj}$ and $\beta_{nj0}$. From (5.22) and (5.23) we obtain

(5.24) $$\|K_4\|^2 \leq 4 \left\{ \sum_{j=1}^{s_n} p'_{\lambda_n}(|\hat{\beta}_{n0j}|)^2 + C\|\hat{\beta}_{n1} - \beta_{n0}\|^2 \right\}^2$$
$$\leq 4 \left\{ p_n a_n^2 + O_p\left(\frac{p_n}{n}\right) \right\}^2 = O_p\left(\frac{p_n^2}{n^2}\right) = o_p(1).$$

Finally, we consider $K_3$. It is easy to see that $K_3$ can be decomposed as the sum of $K_5$ and $K_6$, where

$$K_5 \triangleq \left\{ \frac{1}{n} \sum_{i=1}^n \frac{\partial L_{ni}(\hat{\beta}_{n1})}{\partial \beta_j} \frac{\partial L_{ni}(\hat{\beta}_{n1})}{\partial \beta_k} - \frac{1}{n} \sum_{i=1}^n \frac{\partial L_{ni}(\beta_{n01})}{\partial \beta_j} \frac{\partial L_{ni}(\beta_{n01})}{\partial \beta_k} \right\},$$

$$K_6 \triangleq \left\{ \frac{1}{n} \sum_{i=1}^n \frac{\partial L_{ni}(\beta_{n01})}{\partial \beta_j} \frac{\partial L_{ni}(\beta_{n01})}{\partial \beta_k} \right\} - I_n(\beta_{n01}).$$

As before, following Lemma 8 in the Appendix, it is easy to demonstrate that

(5.25) $$\|K_6\| = o_p(1).$$

In the Appendix we show that

(5.26) $$\|K_5\| = o_p(1).$$

By (5.24)–(5.26) we have shown that $\|\mathcal{B}_n - \mathcal{B}\| = o_p(1)$ and

(5.27) $$\lambda_i(\mathcal{B}_n - \mathcal{B}) = o_p(1), \qquad i = 1, \ldots, s_n.$$

It follows from (5.21) and (5.27) that

$$\lambda_i(\hat{\Sigma}_n - \Sigma_n) = o_p(1), \qquad i = 1, \ldots, s_n.$$

This completes the proof for the consistency of the sandwich formula. □

Let $B_n$ be an $(s_n - q) \times s_n$ matrix which satisfies $B_n B_n^T = I_{s_n - q}$ and $A_n B_n^T = 0$. As $\beta_{n1}$ is in the orthogonal complement to the linear space that is spanned by rows of $A_n$ under the null hypothesis $H_0$, it follows that

$$\beta_n = B_n^T \gamma,$$

where $\gamma$ is an $(s_n - q) \times 1$ vector. Then, under $H_0$ the penalized likelihood estimator is also the local maximizer $\hat{\gamma}_n$ of the problem

$$Q_n(B_n^T \hat{\gamma}_n) = \max_{\gamma_n} Q_n(B_n^T \gamma_n).$$

To prove Theorem 4 we need the following two lemmas, the proofs of which are given in the Appendix.

LEMMA 6. *Under the conditions of Theorem 4 and the null hypothesis $H_0$, we have*

$$\hat{\beta}_{n1} - \beta_{n01} = \frac{1}{n} I_n(\beta_{n01})^{-1} \nabla L_n(\beta_{n01}) + o_p(n^{-1/2}),$$

$$B_n^T(\hat{\gamma}_n - \gamma_{n0}) = \frac{1}{n} B_n^T \{B_n I_n(\beta_{n01}) B_n^T\}^{-1} B_n^T \nabla L_n(\beta_{n01}) + o_p(n^{-1/2}).$$

LEMMA 7. *Under the conditions of Theorem 4 and the null hypothesis $H_0$, we have*

$$\begin{aligned}(5.28)\quad & Q_n(\hat{\beta}_{n1}) - Q_n(B_n^T \hat{\gamma}_n) \\ & = \frac{n}{2}(\hat{\beta}_{n1} - B_n^T \hat{\gamma}_n)^T I_n(\beta_{n01})(\hat{\beta}_{n1} - B_n^T \hat{\gamma}_n) + o_p(1).\end{aligned}$$

PROOF OF THEOREM 4. Let $\Theta_n = I_n(\beta_{n01})$ and $\Phi_n = \frac{1}{n} \nabla L_n(\beta_{n01})$. By Lemma 6 we have

$$\begin{aligned}(5.29)\quad \hat{\beta}_{n1} - B_n^T \hat{\gamma}_n &= \Theta_n^{-1/2} \{I_n - \Theta_n^{1/2} B_n^T (B_n \Theta_n B_n^T)^{-1} B_n \Theta_n^{1/2}\} \Theta_n^{-1/2} \Phi_n \\ &\quad + o_p(n^{-1/2}).\end{aligned}$$

It is easy to see that $I_n - \Theta_n^{1/2} B_n^T (B_n \Theta_n B_n^T)^{-1} B_n \Theta_n^{1/2}$ is an idempotent matrix with rank $q$. Hence, by a standard argument and condition (F),

$$\hat{\beta}_{n1} - B_n^T \hat{\gamma}_n = O_p\left(\sqrt{\frac{q}{n}}\right).$$

Substituting (5.29) into (5.28), we obtain

$$\begin{aligned}Q_n(\hat{\beta}_{n1}) - Q_n(B_n^T \hat{\gamma}_n) \\ = \frac{n}{2} \Phi_n^T \Theta_n^{-1/2} \{I_n - \Theta_n^{1/2} B_n^T (B_n \Theta_n B_n^T)^{-1} B_n \Theta_n^{1/2}\} \Theta_n^{-1/2} \Phi_n + o_p(1).\end{aligned}$$

By the property of the idempotent matrix, $I_n - \Theta_n^{1/2} B_n^T (B_n \Theta_n B_n^T)^{-1} B_n \Theta_n^{1/2}$ can be written as the product form $D_n^T D_n$, where $D_n$ is a $q \times s_n$ matrix that satisfies $D_n D_n^T = I_q$. As in the proof of Theorem 2, we have shown that $\sqrt{n} D_n \Theta_n^{-1/2} \Phi_n$ has an asymptotic multivariate normal distribution, that is,

$$\sqrt{n} D_n \Theta_n^{-1/2} \Phi_n \xrightarrow{\mathcal{D}} \mathcal{N}(0, I_q).$$

Finally, we have

$$\begin{aligned}2\{Q_n(\hat{\beta}_{n1}) - Q_n(B_n^T \hat{\gamma}_n)\} \\ = n(D_n \Theta_n^{-1/2} \Phi_n)^T (D_n \Theta_n^{-1/2} \Phi_n) + o_p(1) \xrightarrow{\mathcal{D}} \chi_q. \quad \square\end{aligned}$$

# APPENDIX

LEMMA 8. *Under the conditions of Theorem* 1, *we have*

(A.1) $$\left\| \frac{1}{n} \nabla^2 L_n(\beta_{n0}) + I_n(\beta_{n0}) \right\| = o_p\left(\frac{1}{p_n}\right)$$

*and*

(A.2) $$\left\| \left\{ \frac{1}{n} \sum_{i=1}^{n} \frac{\partial L_{ni}(\beta_{n01})}{\partial \beta_j} \frac{\partial L_{ni}(\beta_{n01})}{\partial \beta_k} \right\} - I_n(\beta_{n0}) \right\| = o_p\left(\frac{1}{p_n}\right).$$

PROOF. For any $\varepsilon$, by Chebyshev's inequality,

$$P\left(\left\| \frac{1}{n} \nabla^2 L_n(\beta_{n0}) + I_n(\beta_{n0}) \right\| \geq \frac{\varepsilon}{p_n}\right)$$
$$\leq \frac{p_n^2}{n^2 \varepsilon^2} E \sum_{i,j=1}^{p_n} \left\{ \frac{\partial L_n(\beta_{n0})}{\partial \beta_{ni} \beta_{nj}} - E \frac{\partial L_n(\beta_{n0})}{\partial \beta_{ni} \beta_{nj}} \right\}^2$$
$$= \frac{p_n^4}{n} = o(1).$$

Hence (A.1) follows. Similarly, we can prove (A.2). □

LEMMA 9. *Under the conditions of Theorem* 2, *we have*

$$\left\| \frac{\nabla^2 L_n(\hat{\beta}_{n1})}{n} - \frac{\nabla^2 L_n(\beta_{n01})}{n} \right\| = o_p\left(\frac{1}{\sqrt{p_n}}\right).$$

PROOF. First we expand the left-hand side of the equation above to the third order,

$$\left\| \frac{\nabla^2 L_n(\hat{\beta}_{n1})}{n} - \frac{\nabla^2 L_n(\beta_{n01})}{n} \right\|^2$$
$$= \frac{1}{n^2} \sum_{i,j=1}^{s_n} \left\{ \frac{\partial L_n(\hat{\beta}_{n1})}{\partial \beta_i \, \partial \beta_j} - \frac{\partial L_n(\beta_{n01})}{\partial \beta_i \, \partial \beta_j} \right\}^2$$
$$= \frac{1}{n^2} \sum_{i,j=1}^{s_n} \left\{ \sum_{k=1}^{s_n} \frac{\partial L_n(\beta_{n1}^*)}{\partial \beta_i \, \partial \beta_j \, \partial \beta_k} (\hat{\beta}_{nk} - \beta_{n0k}) \right\}^2.$$

Then by condition (G) and the Cauchy–Schwarz inequality,

$$\frac{1}{n^2} \sum_{i,j=1}^{s_n} \left\{ \sum_{k=1}^{s_n} \frac{\partial L_n(\beta_{n1}^*)}{\partial \beta_i \, \partial \beta_j \, \partial \beta_k} (\hat{\beta}_{nk} - \beta_{n0k}) \right\}^2$$

$$\leq \frac{1}{n^2} \sum_{i,j=1}^{s_n} \sum_{k=1}^{s_n} \left\{ \frac{\partial L_n(\beta_{n1}^*)}{\partial \beta_i \, \partial \beta_j \, \partial \beta_k} \right\}^2 \|\hat{\beta}_{n1} - \beta_{n01}\|^2$$

$$= \frac{1}{n^2} O_p\left(\frac{p_n}{n}\right) \sum_{i,j,k}^{p_n} \left\{ \sum_{l=1}^{n} M_{nijk}(V_{nl}) \right\}^2$$

$$= \frac{1}{n^2} O_p\left(\frac{p_n}{n}\right) O_p(p_n^3 n^2) = o_p\left(\frac{1}{p_n}\right). \qquad \Box$$

PROOF OF (5.26). According to Taylor's expansion, we have

$$\frac{\partial L_{ni}(\hat{\beta}_{n1})}{\partial \beta_j} = \frac{\partial L_{ni}(\beta_{n01})}{\partial \beta_j} + \nabla^T \frac{\partial L_{ni}(\beta_{n01})}{\partial \beta_j} (\hat{\beta}_{n1} - \beta_{n0})$$

$$+ (\hat{\beta}_{n1} - \beta_{n0})^T \nabla^2 \frac{\partial L_{ni}(\beta_{n1}^*)}{\partial \beta_j} (\hat{\beta}_{n1} - \beta_{n0})$$

$$\hat{=} a_{ij} + b_{ij} + c_{ij}.$$

The matrix $K_5$ can then be expressed as a sum of the following form:

$$K_5 = \frac{1}{n}\left(\sum_{i=1}^{n} a_{ij}b_{ik}\right) + \frac{1}{n}\left(\sum_{i=1}^{n} a_{ij}c_{ik}\right) + \frac{1}{n}\left(\sum_{i=1}^{n} b_{ij}a_{ik}\right) + \frac{1}{n}\left(\sum_{i=1}^{n} c_{ij}a_{ik}\right)$$

$$+ \frac{1}{n}\left(\sum_{i=1}^{n} b_{ij}b_{ik}\right) + \frac{1}{n}\left(\sum_{i=1}^{n} b_{ij}c_{ik}\right) + \frac{1}{n}\left(\sum_{i=1}^{n} c_{ij}b_{ik}\right) + \frac{1}{n}\left(\sum_{i=1}^{n} c_{ij}c_{ik}\right)$$

$$\hat{=} X_1 + X_2 + X_3 + X_4 + X_5 + X_6 + X_7 + X_8.$$

Considering a matrix of the form $n^{-1}(\sum_{i=1}^{n} x_{ij}y_{ik}) \hat{=} F$, we have

$$\|F\|^2 = \frac{1}{n^2} \sum_{j,k=1}^{s_n} \left(\sum_{i=1}^{n} x_{ij}y_{ik}\right)^2 \leq \frac{1}{n^2} \sum_{j,k=1}^{s_n} \left(\sum_{i=1}^{n} x_{ij}^2\right)\left(\sum_{i=1}^{n} y_{ik}^2\right)$$

$$= \frac{1}{n^2}\left(\sum_{i=1}^{n} \sum_{j=1}^{s_n} x_{ij}^2\right)\left(\sum_{i=1}^{n} \sum_{k=1}^{s_n} y_{ik}^2\right).$$

Thus, the order of $\|X_i\|$ can be determined from those of $\sum_{i=1}^{n}\sum_{j=1}^{s_n} a_{ij}^2$, $\sum_{i=1}^{n}\sum_{j=1}^{s_n} b_{ij}^2$ and $\sum_{i=1}^{n}\sum_{j=1}^{s_n} c_{ij}^2$.

Because of condition (F), for any $i$ and $j$ $Ea_{ij}^2 \leq C$ and

$$E\left\{\frac{\partial L_{ni}(\beta_{n01})}{\partial \beta_j \, \partial \beta_k}\right\}^2 \leq C \qquad \text{for any } n, j \text{ and } k,$$

we obtain

(A.3) $$\sum_{i=1}^{n} \sum_{j=1}^{s_n} a_{ij}^2 = O_p(np_n)$$

and

$$\sum_{i=1}^{n}\sum_{j=1}^{s_n} b_{ij}^2 \leq \sum_{i=1}^{n}\sum_{j=1}^{s_n}\sum_{k=1}^{s_n}\left\{\frac{\partial L_{ni}(\beta_{n01})}{\partial\beta_j\,\partial\beta_k}\right\}^2 \|\hat{\beta}_{n1} - \beta_{n01}\|^2$$

(A.4)
$$= O_p(np_n^2) O_p\left(\frac{p_n}{n}\right) = O_p(p_n^3).$$

By condition (G) and using the Cauchy–Schwarz inequality, we show that

$$\sum_{i=1}^{n}\sum_{j=1}^{s_n} c_{ij}^2 \leq \sum_{i=1}^{n}\sum_{j=1}^{s_n}\sum_{k=1}^{s_n}\sum_{l=1}^{s_n}\left\{\frac{\partial L_{ni}(\beta_{n1}^*)}{\partial\beta_j\,\partial\beta_k\,\partial\beta_l}\right\}^2 \|\hat{\beta}_{n1} - \beta_{n01}\|^4$$

(A.5)
$$= O_p(np_n^3) O_p\left(\frac{p_n^2}{n^2}\right) = O_p\left(\frac{p_n^5}{n}\right).$$

From (A.3)–(A.5) we have

$$\|K_5\|^2 \leq 8(\|X_1\|^2 + \cdots + \|X_8\|^2)$$
$$\leq 8\frac{1}{n^2}\bigg\{O_p(np_n p_n^3) + O_p\bigg(np_n \frac{p_n^5}{n}\bigg)$$
$$+ O_p\bigg(p_n^3 \frac{p_n^5}{n}\bigg) + O_p(p_n^3 p_n^3) + O_p\bigg(\frac{p_n^5}{n}\frac{p_n^5}{n}\bigg)\bigg\}$$
$$= O_p\bigg(\frac{p_n^4}{n}\bigg) = o_p(1).$$

This completes the proof. □

PROOF OF LEMMA 6. We need only prove the second equation. The first equation can be shown in the same manner. Following the steps of the proof of Theorem 2, it follows that under $H_0$,

$$B_n(I_n(\beta_{n01}) + \Sigma_{\lambda_n})B_n^T(\hat{\gamma}_n - \gamma_{n0}) - B_n\mathbf{b}_n = \frac{1}{n}B_n\nabla L_n(\beta_{n01}) + o_p(n^{-1/2}).$$

By the conditions $a_n = o_p(1/\sqrt{np_n})$ and $B_n B_n^T = I_{s_n - q}$, we have

$$\|B_n\mathbf{b}_n\| \leq \|\mathbf{b}_n\| \leq \sqrt{p_n}a_n = o_p(n^{-1/2}).$$

On the other hand, since $b_n = o_p(1/\sqrt{p_n})$, we obtain

$$\|B_n\Sigma_{\lambda_n}B_n^T(\hat{\gamma}_n - \gamma_{n0})\| \leq \|\hat{\gamma}_n - \gamma_{n0}\|b_n = o_p\left(\frac{1}{\sqrt{p_n}}\right)O_p\left(\sqrt{\frac{p_n}{n}}\right) = o_p\left(\frac{1}{\sqrt{n}}\right).$$

Hence, it follows that

$$B_n I_n(\beta_{n01})B_n^T(\hat{\gamma}_n - \gamma_{n0}) = \frac{1}{n}B_n\nabla L_n(\beta_{n01}) + o_p(n^{-1/2}).$$

As $\lambda_i(B_n I_n(\beta_{n01}) B_n^T)$ is uniformly bounded away from 0 and infinity, we have

$$B_n^T(\hat{\gamma}_n - \gamma_{n0}) = B_n^T \{B_n I_n(\beta_{n01}) B_n^T\}^{-1} B_n \nabla L_n(\beta_{n01}) + o_p(n^{-1/2}).$$

This completes the proof. □

PROOF OF LEMMA 7. A Taylor's expansion of $Q_n(\hat{\beta}_{n1}) - Q_n(B_n^T \hat{\gamma}_n)$ at the point $\hat{\beta}_{n1}$ yields

$$Q_n(\hat{\beta}_{n1}) - Q_n(B_n^T \hat{\gamma}_n) = T_1 + T_2 + T_3 + T_4,$$

where

$$T_1 = \nabla^T Q_n(\hat{\beta}_{n1})(\hat{\beta}_{n1} - B_n^T \hat{\gamma}),$$
$$T_2 = -\tfrac{1}{2}(\hat{\beta}_{n1} - B_n^T \hat{\gamma}_n)^T \nabla^2 L_n(\hat{\beta}_{n1})(\hat{\beta}_{n1} - B_n^T \hat{\gamma}_n),$$
$$T_3 = \tfrac{1}{6}\nabla^T \{(\hat{\beta}_{n1} - B_n^T \hat{\gamma}_n)^T \nabla^2 L_n(\beta_{n1}^*)(\hat{\beta}_{n1} - B_n^T \hat{\gamma}_n)\}(\hat{\beta}_{n1} - B_n^T \hat{\gamma}_n),$$
$$T_4 = \tfrac{1}{2}(\hat{\beta}_{n1} - B_n^T \hat{\gamma}_n)^T \nabla^2 P_{\lambda_n}(\hat{\beta}_{n1}) \{I + o(I)\}(\hat{\beta}_{n1} - B_n^T \hat{\gamma}_n).$$

Note that $T_1 = 0$ as $\nabla^T Q(\hat{\beta}_{n1}) = 0$. By Lemma 6 and (5.29) it follows that

$$\hat{\beta}_{n1} - B_n^T \hat{\gamma}_n = O_p\left(\sqrt{\frac{q}{n}}\right).$$

By the conditions $b_n = o_p(1/\sqrt{p_n})$ and $q < p_n$, following the proof of $I_3$ in Theorem 1, we have

$$T_3 = O_p(n p_n^{3/2} n^{-3/2} q^{3/2}) = o_p(1)$$

and

$$T_4 \le n b_n \|\hat{\beta}_{n1} - B_n^T \hat{\gamma}_n\|^2 = n\, o_p\left(\frac{1}{\sqrt{p_n}}\right) O_p\left(\frac{q}{n}\right) = o_p(1).$$

Thus,

(A.6) $$Q_n(\hat{\beta}_{n1}) - Q_n(B_n^T \hat{\gamma}_n) = T_2 + o_p(1).$$

It follows from Lemmas 8 and 9 that

$$\left\|\frac{1}{n} \nabla^2 L_n(\hat{\beta}_{n1}) + I_n(\beta_{n01})\right\| = o_p\left(\frac{1}{\sqrt{p_n}}\right).$$

Hence, we have

$$\frac{1}{2}(\hat{\beta}_{n1} - B_n^T \hat{\gamma}_n)^T \{\nabla^2 L_n(\hat{\beta}_{n1}) + n I_n(\beta_{n01})\}(\hat{\beta}_{n1} - B_n^T \hat{\gamma}_n)$$
$$\le o_p\left(n \frac{1}{\sqrt{p_n}}\right) O_p\left(\frac{q}{n}\right) = o_p(1).$$

The combination of (A.6) and (A.7) yields (5.28). □

**Acknowledgments.** The authors gratefully acknowledge the helpful comments of the Associate Editor and referees.

DEPARTMENT OF OPERATIONS RESEARCH
AND FINANCIAL ENGINEERING
PRINCETON UNIVERSITY
PRINCETON, NEW JERSEY 08544
USA
E-MAIL: jqfan@princeton.edu

DEPARTMENT OF STATISTICS
THE CHINESE UNIVERSITY OF HONG KONG
SHATIN, HONG KONG
E-MAIL: h_peng@sparc20c.sta.cuhk.edu.hk